\documentclass[3p,times]{elsarticle}
\usepackage{amsfonts}

\journal{Journal of Mathematical Analysis and Applications (be published in 2016)}

\usepackage{amssymb}
\usepackage{amsmath}
\usepackage{mathrsfs}
\numberwithin{equation}{section}
\newtheorem{thm}{Theorem}[section]
\newtheorem{lem}{Lemma}[section]
\newdefinition{rem}{Remark}[section]
\newtheorem{prop}{Proposition}[section]
\newtheorem{cor}{Corollary}[section]
\newdefinition{defi}{Definition}[section]
\newproof{pf}{Proof}
\usepackage[figuresright]{rotating}

\begin{document}

\begin{frontmatter}



%% e.g. \dochead{17th International Conference on Dynamical Processes in Excited States of Solids}

\title{Existence and construction of quasi-stationary distributions for one-dimensional diffusions}


\author[label1]{Hanjun Zhang}
\author[label1]{Guoman He\corref{cor1}}
\ead{hgm0164@163.com}
\cortext[cor1]{Corresponding author}
\address[label1]{School of Mathematics and Computational Science, Xiangtan University, Hunan 411105, P. R. China.}

\begin{abstract}
In this paper, we study quasi-stationary distributions (QSDs) for one-dimensional diffusions killed at 0, when 0 is a regular boundary and $+\infty$ is a natural boundary. More precisely, we not only give a necessary and sufficient condition for the existence of a QSD, but we also construct all QSDs for the one-dimensional diffusions. Moreover, we give a sufficient condition for $R$-positivity of the process killed
at the origin. This condition is only based on the drift, which is easy to check.
\end{abstract}

\begin{keyword}
One-dimensional diffusion; Quasi-stationary distribution; $R$-positivity; Natural boundary

\MSC primary 60J60; secondary 60J70; 37A30
\end{keyword}

\end{frontmatter}

% \linenumbers

\section{Introduction}
\label{sect1}
We are interested in the long-term behavior of killed Markov processes. Conditional stationarity, which we call quasi-stationarity, is one of the most interesting topics in this direction. For quasi-stationary distribution (QSD), we know that the study of QSDs is a long standing problem in several areas of probability theory
and a complete understanding of the structure of QSDs seems to be available only in rather special situations such as Markov chains on finite sets or more general processes with compact state space. The main motivation of this work is the existence and construction of QSDs for one-dimensional diffusion $X$ killed at 0, when 0 is a regular boundary and $+\infty$ is a natural boundary. Moreover, we give a sufficient condition in order for the process $X$ killed at 0 to be $R$-positive.
\par
To the best of our knowledge, Mandl \cite{M61} is the first one to study the existence of a QSD for continuous time diffusion process on the half line. If Mandl's conditions are satisfied, the existence of the Yaglom limit and that of a QSD for killed one-dimensional diffusion processes have been proved by various authors (see, e.g., \cite{CMS95, KS12, LM00, MPM98, SE07}). If Mandl's conditions are not satisfied, Cattiaux et~al. who studied the existence and uniqueness of the QSD for one-dimensional diffusions killed at 0 and whose drift is allowed to go to $-\infty$ at 0 and the process is allowed to have an entrance boundary at $+\infty$, have done a pioneering work (see \cite{CCLMMS09}). In this case, under the most general conditions, Littin proves the existence of a unique QSD and of the Yaglom limit in \cite{LC12}, which is closely related to \cite{CCLMMS09}. Although \cite{CMS95, KS12, LM00, MPM98, MM04, SE07} and \cite{CCLMMS09} make the key contributions, the structure of QSDs of killed one-dimensional diffusions has not been completely clarified. This leads us to further study QSDs for one-dimensional diffusions.
\par
Another notion is $R$-positivity, which, in general, is not easy to check, is sufficient to facilitate the straightforward calculation of QSDs for a process from the eigenvectors, eigenmeasures and eigenvalues of its transition rate matrix. The classification of killed one-dimensional diffusions has been studied by Mart\'{\i}nez and San Mart\'{\i}n \cite{MM04}. They gave necessary and sufficient conditions, in terms of the bottom eigenvalue function, for $R$-recurrence and $R$-positivity of one-dimensional diffusions killed at the origin.
\par
In this paper, the main novelty is that we not only give a necessary and sufficient condition for the existence of a QSD, but we also construct all QSDs for one-dimensional diffusion $X$ killed at 0, when 0 is a regular boundary and $+\infty$ is a natural boundary. Moreover, compared with \cite{MM04}, we give an explicit criterion for the process $X$ killed at 0 is $R$-positive.
\par
The remainder of this paper is organized as follows. In Section \ref{sect2} we present some preliminaries that will be needed in the sequel. In Section \ref{sect3} we characterize all QSDs for one-dimensional diffusion $X$ killed at 0, when 0 is a regular boundary and $+\infty$ is a natural boundary. In Section \ref{sect4} we mainly show under what direct conditions on the drift the process is $R$-positive. We conclude in Section \ref{sect5} with some examples.
\section{Preliminaries}
\label{sect2}
We consider the generator ${L}u={1\over2}\partial_{xx}u-q\partial_xu$. Denote by $X$ the diffusion whose infinitesimal generator is $L$, or in other words the solution of the stochastic differential equation (SDE)
\begin{equation}
\label{2.1}
dX_t=dB_t-q(X_t)dt,~~~~~~~ X_0=x>0,
\end{equation}
where $(B_t;t\geq0)$ is a standard one-dimensional Brownian motion and $q\in C^1([0,\infty))$. Thus, $-q$ is the drift of $X$. Observe that, under the condition $q\in C^1([0,\infty))$, $\int_{0}^{d}e^{Q(y)}dy<\infty$ and $\int_{0}^{d}e^{-Q(y)}dy<\infty$ for some (and, therefore, for all) $d>0$, which is equivalent to saying that the boundary point 0 is regular in the sense of Feller, where $Q(y)=\int_{0}^{y}2q(x)dx$.
\par
Let $\mathbb{P}_x$ and $\mathbb{E}_x$ stand for the probability and the expectation, respectively, associated with $X$ when initiated from $x$.
For any probability measure $\nu$ on $(0,\infty)$, we denote by $\mathbb{P}_\nu$ the probability associated with the process $X$ initially
distributed with respect to $\nu$. Let $\tau_a:=\inf\{t>0:X_t=a\}$ be the hitting time of $a$. We are mainly interested in the case $a=0$ and we denote $\tau=\tau_0$. As usual $X^\tau$ corresponds to $X$ killed at 0.
\par
Associated with $q$, we consider the function
\begin{equation}
\label{2.2}
\Lambda(x)=\int_{0}^{x}e^{Q(y)}dy.
\end{equation}
Notice that
$\Lambda$ is the scale  function for $X$. It satisfies $L\Lambda\equiv0,~\Lambda({0})=0,~\Lambda'({0})=1$. It will be useful to introduce the following measure defined on $(0,\infty)$:
\begin{equation}
\label{2.3}
\mu(dy):=e^{-Q(y)}dy.
\end{equation}
Notice that $\mu$ is the speed measure for $X$.
 \par
Let ${L}^*={1\over2}\partial_{xx}+\partial_x(q \cdot)$
be the formal adjoint operator of $L$. We denote by $\varphi_\lambda$ the solution of
\begin{equation}
\label{2.4}
{L}^*\varphi_\lambda=-\lambda\varphi_\lambda,~~\varphi_\lambda({0})=0,~~\varphi_\lambda'({0})=1,
\end{equation}
and by $\eta_\lambda$ the solution of
\begin{equation}
\label{2.5}
{L}\eta_\lambda=-\lambda\eta_\lambda,~~\eta_\lambda({0})=0,~~\eta_\lambda'({0})=1.
\end{equation}
A direct computation shows that
\begin{equation}
\label{2.6}
\varphi_\lambda=e^{-Q}\eta_\lambda.
\end{equation}
\par
Let $\lambda_c$ be the smallest point of increase of $\varrho(\lambda)$, where $\varrho(\lambda)$ denotes the spectral measure of the operator ${L}^*$. We will assume $\varrho(\lambda)$ is left-continuous (see \cite[Chapter 9]{CL55}).
In \cite[Lemma 2]{M61} it was shown that
$$\lambda_c=\sup\{\lambda:\varphi_\lambda(\cdot) ~\mathrm{does~ not~ change ~sign}\}.$$
\par
For most of the results in this paper we will use the following hypothesis $(\mathrm{H})$, that is,
\begin{defi}
\label{defi1}
We say that hypothesis $(\mathrm{H})$ holds if the following explicit conditions on $q$, all together, are satisfied:
\par $(\mathrm{H1})$ $\Lambda(\infty)=\infty$.
\par $(\mathrm{H2})$ $S=\int_{0}^{\infty}e^{Q(y)}\left(\int_y^{\infty}e^{-Q(z)}dz\right)dy=\infty$.
\end{defi}
\par If (H1) holds, then it is equivalent to $\mathbb{P}_x(\tau<\infty)=1$, for all $x>0$ (see, e.g., \cite[Chapter VI, Theorem 3.2]{IW89}). So that, if (H1) and (H2) are satisfied, then $+\infty$ is a natural boundary according to Feller's classification (see \cite[Chapter 15]{KT81}).
\par

\section{Existence and construction of quasi-stationary distributions}
\label{sect3}
In this section, we study the standard quasi-stationary distributions of a one-dimensional diffusion $X$ killed at 0, when 0 is a regular boundary and $+\infty$ is a natural boundary, a typical problem for absorbing Markov processes. More formally, the following definition captures the main object of interest of this work.
\par
\begin{defi}
\label{defi5}
We say that a probability measure $\nu$ supported on $(0,\infty)$ is a quasi-stationary distribution $\mathrm{(QSD)}$, if for all $t\geq0$ and any Borel subset $A$ of $(0,\infty)$,
\begin{equation}
\label{3.1}
 \mathbb{P}_\nu(X_t\in A|\tau>t)=\nu(A).
 \end{equation}
\end{defi}
\par
It is well known that a basic and useful property is that the time of killing is exponentially distributed when starting from a QSD:
\begin{prop}
\label{prop1}
Assume that $\pi$ is a $\mathrm{QSD}$ for the process $X$. Then there exists $\lambda>0$ such that, for all $t>0$,
\begin{equation}
\label{3.2}
\mathbb{P}_\pi(\tau>t)=e^{-\lambda t}.
\end{equation}
\end{prop}
\par
The following theorem is one of our main results.
\begin{thm}
\label{Theorem 1}
There exists a quasi-stationary distribution for one-dimensional diffusion $X$ satisfying $(\ref{2.1})$ if and only if $(\mathrm{H1})$ is satisfied and the following condition$:$
\begin{equation}
\label{3.3}
\delta\equiv\sup_{x>0}\int_{0}^{x}e^{Q(y)}dy\int_x^{\infty}2e^{-Q(y)}dy<\infty
\end{equation}
holds. Moreover, if $(\mathrm{H2})$ holds, then for any $0<\lambda\leq\lambda_c, d\nu_\lambda=2\lambda\eta_\lambda d\mu$ is a quasi-stationary distribution and all of quasi-stationary distributions for $X$ only have this form, where $\mu$ and $\eta_\lambda$ are defined by $(\ref{2.3})$ and $(\ref{2.5})$ respectively.
\end{thm}
\begin{rem}
Assuming that the measure $\mu$ is finite, Chen \cite{MFC00} gave an explicit bounds of the first Dirichlet eigenvalue for the generator $L$ on the half line $\mathbb{R}^+:=[0,\infty)$ in terms of $\delta$. Assuming that absorption is certain, i.e. (H1) holds, Pinsky \cite{PRG09} also obtained the same sharp estimate for the first Dirichlet eigenvalue of $L$ on $\mathbb{R}^+$. Now we firstly use the condition (\ref{3.3}) to study the existence of QSDs for $X$.
\end{rem}
\par A relevant quantity in our study is the exponential decay for the absorption probability
$$\zeta=-\lim_{t\rightarrow\infty}\frac{1}{t}\log\mathbb{P}_x(\tau>t).$$
 \par In \cite[Theorem A]{CMS95} it was shown that $\zeta$ exists and is independent of $x>0$. Noticing that, the following three lemmas (Lemmas \ref{Lem 3.1}--\ref{Lem 3.3}) play a key role in this paper, which had been proved in \cite{MM01}.
\begin{lem}
\label{Lem 3.1}
Assume $(\mathrm{H2})$ holds. The following properties are equivalent$:$\\
$~~~~~~~(\mathrm{i})$ $\zeta>0$$;$\\
$~~~~~~(\mathrm{ii})$ $\int_{0}^{\infty}\varphi_{\lambda_c}(y)dy<\infty$ and $\int_{0}^{\infty}e^{-Q(y)}dy<\infty$$;$\\
$~~~~~~(\mathrm{iii})$ $\int_{0}^{\infty}\varphi_{\lambda_c}(y)dy<\infty$ and $\Lambda(\infty)=\infty$$;$\\
$~~~~~~(\mathrm{iv})$ $\lambda_c>0~~and~~\Lambda(\infty)=\infty$$;$\\
$~~~~~~(\mathrm{v})$ $\exists~\lambda>0$ such that $\eta_\lambda$ is increasing.
\end{lem}
\begin{lem}
\label{Lem 3.2}
Assume $(\mathrm{H})$ holds. Then $\zeta=\lambda_c$.
\end{lem}
\begin{lem}
\label{Lem 3.3}
Assume $(\mathrm{H})$ holds. The following statements are equivalent for $\lambda\in(0,\lambda_c]$$:$\\
$~~~~~~(\mathrm{i})$ $\eta_\lambda$ $($or equivalent $\varphi_\lambda$$)$ is positive$;$\\
$~~~~~(\mathrm{ii})$ $\eta_\lambda$ is strictly increasing$;$\\
$~~~~(\mathrm{iii})$ $\varphi_\lambda$ is  strictly positive and integrable.\\
Moreover, if any of these conditions holds, then
\begin{equation}
\label{3.4}
\lim_{y\rightarrow\infty}\frac{\eta_\lambda(y)}{\Lambda(y)}=0~~~~and~~~~\int_{0}^{\infty}\varphi_{\lambda}(x)dx=\frac{1}{2\lambda}.
\end{equation}
\end{lem}
\par
To yield our result more conveniently, we introduce the following lemma.
\begin{lem}
\label{Lem 3.4}
Assume $(\mathrm{H})$ holds. Then $\lambda_c>0$ if and only if $\delta<\infty$, where $\delta$ is defined by $(\ref{3.3})$.
\end{lem}
\begin{pf}
Assume $(\mathrm{H})$ holds. If $\lambda_c>0$, we know from Lemma \ref{Lem 3.1} that the measure $\mu$ is finite.
According to \cite[Theorem 1.1]{MFC00}, we know that
$$(4\delta)^{-1}\leq\lambda_c\leq(\delta)^{-1}.$$
Thus we deduce $\delta<\infty$.
\par
Conversely, we have assumed that 0 is a regular boundary for the process $X$. So if $\delta<\infty$, we can deduce
the measure $\mu$ is finite. By using \cite[Theorem 1.1]{MFC00} again, we deduce $\lambda_c>0$.
\qed
\end{pf}
\par From $(\ref{2.3})$ and $(\ref{2.6})$, we can define $$d\nu_\lambda=2\lambda\varphi_\lambda(y)dy=2\lambda\eta_\lambda(y)e^{-Q(y)}dy
=2\lambda\eta_\lambda(y)\mu(dy)=2\lambda\eta_\lambda d\mu.$$
\par We may now study the existence of QSD and construct all QSDs under the condition $(\mathrm{H})$ is satisfied, which is equivalent to the property that $+\infty$ is a natural boundary for $X$ in the sense of Feller. The result is presented in the following proposition.
\begin{prop}
\label{pro 3.1}
Assume both $(\mathrm{H1})$ and $(\ref{3.3})$ hold. Then for any $\lambda\in(0,\lambda_c]$,
$$d\nu_\lambda=2\lambda\eta_\lambda d\mu$$
is a quasi-stationary distribution if and only if the following two conditions are satisfied$:$\\
$~~~~~~(\mathrm{i})$ $\int_{0}^{\infty}d\nu_\lambda=1;$\\
$~~~~~(\mathrm{ii})$ ${L}^*\nu_\lambda=-\lambda\nu_\lambda.$
\end{prop}
\begin{pf}
Thanks to the equality (\ref{3.4}), we know that $\nu_\lambda$ is a probability measure, i.e. $\int_{0}^{\infty}d\nu_\lambda=\int_{0}^{\infty}2\lambda\eta_\lambda d\mu=1$. Hence, the condition (i) is satisfied.
\par
Moreover, we know from the equality (\ref{2.4}) that
$${L}^*\nu_\lambda={L}^*2\lambda\varphi_\lambda=-2\lambda^2\varphi_\lambda=-\lambda\nu_\lambda.$$
Therefore, the condition (ii) is satisfied. Next, we will prove that $\nu_\lambda$ is a QSD.
\par
According to \cite{MM01}, for $\lambda\in(0,\lambda_c]$, $\eta_\lambda$ as a solution of the equation $(\ref{2.5})$ is well-defined. Also we know from \cite[Lemma 4]{MM01} that when $\zeta>0$ any of the solutions, $\eta_\lambda$, $\lambda\in(0,\zeta]$, satisfies the semigroup property
\begin{equation}
\label{3.5}
P_t\eta_\lambda(x)=e^{-\lambda t}\eta_\lambda(x)~~~~~~~~\mathrm{for~all}~~x>0, t\geq0.
\end{equation}
Note that $\zeta>0$ and $\lambda_c=\zeta$ can be guaranteed under the conditions $(\mathrm{H1})$ and $(\ref{3.3})$.
\par
Since both $(\mathrm{H1})$ and (\ref{3.3}) hold, thus we know from Lemmas \ref{Lem 3.1}--\ref{Lem 3.4} that
$$\int_{0}^{\infty}\varphi_{\lambda}(y)dy=\int_{0}^{\infty}\eta_\lambda(y)e^{-Q(y)}dy=\int_{0}^{\infty}\eta_\lambda(y)\mu(dy)=\frac{1}{2\lambda}<\infty.$$
Then we have $\eta_\lambda\in\mathbb{L}^1(\mu)$. Thanks to the symmetry of the semigroup,
for all $f\in\mathbb{L}^2(\mu)$ we have
\begin{equation}
\label{3.6}
\int P_tf\eta_\lambda d\mu=\int fP_t\eta_\lambda d\mu=e^{-\lambda t}\int f\eta_\lambda d\mu.
\end{equation}
Since $\eta_\lambda\in\mathbb{L}^1(\mu)$, the equality $(\ref{3.6})$ can extend to all bounded function $f$. In particular, we may use it with $f={\bf1}_{(0,\infty)}$ and with $f={\bf1}_A$, where ${\bf1}_A$ is the indicator function of $A$. Note that
$$\int P_t({\bf1}_{(0,\infty)})2\lambda\eta_\lambda d\mu=\mathbb{P}_{\nu_\lambda}(\tau>t)$$
and
$$\int P_t{\bf1}_A2\lambda\eta_\lambda d\mu=\mathbb{P}_{\nu_\lambda}(X_t\in A,\tau>t),$$
then
\begin{eqnarray*}
\mathbb{P}_{\nu_\lambda}(X_t\in A|\tau>t)&=&{{\mathbb{P}_{\nu_\lambda}(X_t\in A,\tau>t)}\over{{\mathbb{P}_{\nu_\lambda}(\tau>t)}}}={{\int P_t{\bf1}_A2\lambda\eta_\lambda d\mu}\over{\int P_t({\bf1}_{(0,\infty)})2\lambda\eta_\lambda d\mu}}\\
&=&{{\int {\bf1}_AP_t2\lambda\eta_\lambda d\mu}\over{\int ({\bf1}_{(0,\infty)})P_t2\lambda\eta_\lambda d\mu}}={{e^{-\lambda t}\int_A2\lambda\eta_\lambda d\mu}\over{e^{-\lambda t}\int_{0}^{\infty}2\lambda\eta_\lambda d\mu}}\\
&=&\nu_\lambda(A).
\end{eqnarray*}
Thus, we get that $\nu_\lambda$ is a QSD.
\par
Conversely, assume that $\nu_\lambda$ is a QSD. From the definition of QSD, we know that $\nu_\lambda$ is a probability measure, then $\nu_\lambda$ satisfies the condition (i).\par
We only denote $\nu=\nu_\lambda$ here. As defined above, a QSD $\nu$ is a probability measure on $(0,\infty)$ such that for every Borel set $A$ of $(0,\infty)$,
\begin{eqnarray*}
\nu(A)=\frac{\mathbb{P}_\nu(X_t\in A, \tau>t)}{\mathbb{P}_\nu(\tau>t)}&=&\frac{\int P_t({\bf1}_A)(x)\nu(dx)}{\int P_t({\bf1}_{(0,\infty)})(x)\nu(dx)}\\
&=& \frac{P^*_t\nu({\bf1}_A)}{P^*_t\nu({\bf1}_{(0,\infty)})}.
\end{eqnarray*}
where $P^*_t\nu$ is the measure on $(0,\infty)$ defined for $f$ measurable and bounded by
$$P^*_t\nu(f)=\int P_tf(x)\nu(dx).$$
 From the equality (\ref{3.2}), we obtain
$$\int P_t({\bf1}_A)(x)\nu(dx)=P^*_t\nu({\bf1}_A)=e^{-\lambda t}\nu(A).$$
Thus the probability measure $\nu$ is an eigenvector for the operator $P^*_t$ (defined on the signed measure vector space), associated with the eigenvalue $e^{-\lambda t}$. It is easy to show that
$$P^*_t\nu=e^{-\lambda t}\nu\Leftrightarrow\nu P_t=e^{-\lambda t}\nu.$$
Then, it is direct to check that
$${L}^*\nu=-\lambda\nu.$$
Thus $\nu_\lambda$ satisfies the condition (ii). We complete the proof.
\qed
\end{pf}
{\bf Proof of Theorem 3.1.} The theorem follows from Lemmas \ref{Lem 3.1}--\ref{Lem 3.4} and Proposition \ref{pro 3.1}.\qed
\vskip0.2cm
\par
Although the following result also follows from \cite[Lemma 3.3]{KS12} or \cite[Theorem 1]{MM01} due to Lemma \ref{Lem 3.4} we know the condition (\ref{3.3}) is just the condition that $\lambda_c>0$, we give a simple direct proof here.
\begin{cor}
\label{cor 3.1}
If there exists a quasi-stationary distribution for the process $X$, then $\mu(0,\infty)<\infty$.
\end{cor}
\begin{pf}
For any $x\in(0,\infty)$, we have
\begin{equation}
\label{3.7}
\mu(0,\infty)=\int_{0}^{\infty}e^{-Q(z)}dz=\int_{0}^{x}e^{-Q(z)}dz+\int_x^{\infty}e^{-Q(z)}dz.
\end{equation}
Under the assumption, we know from Theorem \ref{Theorem 1} that the equality (\ref{3.3}) holds, then for all $x\in(0,\infty)$, $\int_x^{\infty}e^{-Q(z)}dz<\infty$. Observe that under the condition $q\in C^1([0,\infty))$ we have that $\int_{0}^{x}e^{-Q(z)}dz<\infty$, for all $x\in(0,\infty)$. Thus $\mu(0,\infty)<\infty$ follows immediately.
\qed
\end{pf}

\section{$R$-positivity}
\label{sect4}
In this section, we will show that the one-dimensional diffusion $X$ killed at 0 is $R$-positive. This means that the process $Y$, whose law is the conditional law of $X$ to never hit the origin, is positive recurrent.
\par
A direct computation shows that $\eta_\lambda$ introduced in (\ref{2.5}) satisfies:
\begin{equation}
\label{4.1}
\begin{split}
\eta_\lambda'(x)&=e^{Q(x)}\left(1-2\lambda\int_0^{x}\eta_\lambda(z)e^{-Q(z)}dz\right),\\
\eta_\lambda(x)&=\int_0^{x}e^{Q(y)}\left(1-2\lambda\int_0^{y}\eta_\lambda(z)e^{-Q(z)}dz\right)dy.
\end{split}
\end{equation}
 We know from \cite[Theorem B]{CMS95} that for $x>0$ fixed, the following limit exists and defines a diffusion $Y$:
\begin{eqnarray*}
\lim_{t\rightarrow\infty}\mathbb{P}_x(X_s\in A|\tau>t)&=&e^{\lambda_cs}\mathbb{E}_x\left(\frac{\eta_{\lambda_c}(X_s)}{\eta_{\lambda_c}(x)},X_s\in A,\tau>s\right)\\
      &=&\mathbb{P}_x(Y_s\in A).
\end{eqnarray*}
The diffusion $Y$ satisfies the SDE
\begin{equation}
\label{4.2}
dY_t=dB_t-\phi(Y_t)dt ~~~~~~~~\mathrm{where}~~\phi(y)=q(y)-\frac{\eta'_{\lambda_c}(y)}{\eta_{\lambda_c}(y)},
\end{equation}
and it takes values on $(0,\infty)$. In fact, since its drift is of order $1/x$ for $x$ near 0, so it never reaches 0. \par
The connection between the classification of $Y$ and the $R$-classification of the killed diffusion $X^\tau$ is given in the following definition.
\begin{defi}
\label{defi2}
If the process $Y$ is positive recurrent (resp. recurrent, null recurrent, transient), then the process $X^\tau$ is said to be $R$-positive (resp. $R$-recurrent, $R$-null, $R$-transient).
\end{defi}
\par
We remind the reader of the usual definition of ${\lambda_c}$-positivity of a Markov process here. Let $P_t(x,B), x\in E, B\in \mathscr{B}, t\geq0$ be a Markov transition probability semigroup on a general space $(E, \mathscr{B})$, with $\mathscr{B}$ countably generated, satisfying the simultaneous $\psi$-irreducibility condition. We know from \cite{TT79} that the decay parameter ${\lambda_c}$ of $(P_t)$ exists. Let $(P_t)$ be ${\lambda_c}$-recurrent, that is, for all $x\in E$ and $A\in \mathscr{B}^{+}=\{A\in\mathscr{B}:\psi(A)>0\}, \int e^{\lambda_ct}P_t(x,A)dt=\infty$. The semigroup $(P_t)$ is said to be ${\lambda_c}$-positive, if for every $A\in \mathscr{B}^{+}$
\begin{equation*}
\lim_{t\rightarrow\infty}e^{\lambda_ct}P_t(x,A)>0.
\end{equation*}
We remark that Definition \ref{defi2}
is essentially the same with the usual definition of ${\lambda_c}$-positivity of a Markov process. In fact, we know from \cite[Theorem B]{CMS95} that the transition density $p^Y(t,x,y)$ of $Y$ and the transition density $p^X(t,x,y)$ of $X$ have the following relation:
\begin{equation*}
p^Y(t,x,y)=e^{\lambda_ct}\frac{\eta_{\lambda_c}(y)}{\eta_{\lambda_c}(x)}p^X(t,x,y).
\end{equation*}
Thus, this fact can be clarified easily.
\par
We may now state the following result. In addition, we point out that the following result is stronger than this result: if the bottom of the essential spectrum is strictly bigger than the bottom of the spectrum, then one has $R$-positivity. On this weaker result, we can't directly use it to judge whether a process is $R$-positive or not. We emphasize that our result can provide a direct and explicit condition on the drift such that the process is $R$-positive.
\begin{thm}
\label{thm4.1}
Assume $(\mathrm{H})$ holds. Then $X^\tau$ is $R$-positive if
\begin{equation}
\label{4.3}
\lim\limits_{x\rightarrow\infty}\mu([x,\infty))\int_{0}^{x}e^{Q(y)}dy=0.
\end{equation}
\end{thm}
\begin{pf}
If $(\mathrm{H})$ is satisfied, since $\lim\limits_{x\rightarrow\infty}\mu([x,\infty))\int_{0}^{x}e^{Q(y)}dy=0$ can deduce $\delta=\sup_{x>0}\int_{0}^{x}e^{Q(y)}dy\int_x^{\infty}2e^{-Q(y)}dy<\infty$, then we know from Lemma \ref{Lem 3.4} that $\lambda_c>0$. Further, from Lemma \ref{Lem 3.1}, we obtain $\mu(0,\infty)<\infty$. Next, we will prove that $(\ref{4.3})$ is equivalent to
\begin{equation}
\label{4.4}
\lim\limits_{n\rightarrow\infty}\sup\limits_{r>n}\mu([r,\infty))\int_{n}^{r}e^{Q(x)}dx=0.
\end{equation}
In fact, for any $r>n$, we have
\begin{equation*}
\mu([r,\infty))\int_{n}^{r}e^{Q(x)}dx\leq\mu([r,\infty))\int_{0}^{r}e^{Q(x)}dx,
\end{equation*}
which implies
\begin{equation*}
\lim\limits_{n\rightarrow\infty}\sup\limits_{r>n}\mu([r,\infty))\int_{n}^{r}e^{Q(x)}dx\leq\limsup_{r\rightarrow\infty}\mu([r,\infty))\int_{0}^{r}e^{Q(x)}dx=0.
\end{equation*}
Conversely, for any $n>0$, when $x>n$, we have
\begin{eqnarray*}
\mu([x,\infty))\int_{0}^{x}e^{Q(y)}dy&=&\mu([x,\infty))\int_{0}^{n}e^{Q(y)}dy+\mu([x,\infty))\int_{n}^{x}e^{Q(y)}dy\\
                                     &\leq&\mu([x,\infty))\int_{0}^{n}e^{Q(y)}dy+\sup\limits_{x>n}\mu([x,\infty))\int_{n}^{x}e^{Q(x)}dx.
\end{eqnarray*}
By letting $x\rightarrow\infty$ in the above formula, we have
\begin{equation*}
\lim\limits_{x\rightarrow\infty}\mu([x,\infty))\int_{0}^{x}e^{Q(y)}dy\leq\lim\limits_{n\rightarrow\infty}\sup\limits_{r>n}\mu([r,\infty))\int_{n}^{r}e^{Q(x)}dx=0.
\end{equation*}
Then we prove the equivalence.
\par
We know from \cite[Theorem 1]{PRG09} that (\ref{4.4}) is equivalent to $\sigma_{ess}({L})=\emptyset$, where $\sigma_{ess}({L})$ denotes the essential spectrum of ${L}$. Then we know that $-{L}$ has a purely discrete spectrum $0<\lambda_1<\lambda_2<\cdots$, $\lim_{i\rightarrow\infty}\lambda_i=+\infty$ and there exists an orthonormal basis $\{\eta_i\}_{i=1}^{\infty}$ in $\mathbb{L}^2(\mu)$ such that $-{L}\eta_i=\lambda_i\eta_i$. Here, we remind the reader that $\lambda_1=\lambda_c$.
\par
In order to simplify notation, we shall denote $\eta_1=\eta_{\lambda_c}$. For some $c>0$ fixed, we consider the functions
\begin{equation*}
Q^Y(y)=\int_{c}^{y}2\phi(x)dx=Q(y)-Q(c)-2\log(\eta_1(y)/\eta_1(c))
\end{equation*}
and
\begin{equation*}
\Lambda^Y(y)=\int_{c}^{y}e^{Q^Y(z)}dz=\eta_1^2(c)e^{-Q(c)}\int_{c}^{y}\eta_1^{-2}(z)e^{Q(z)}dz.
\end{equation*}
Because $\eta_1(x)=x+O(x^2)$ for $x$ near 0, we first obtain that $\Lambda^Y(0^+)=-\infty$.
\par The speed measure $m$ of $Y$ is given by
$$m(dx)=\frac{2dx}{(\Lambda^Y(x))'}$$
(see \cite[formula (5.51)]{KS88}). So we obtain
\begin{equation}
\label{4.5}
m(dx)=2\frac{e^{Q(c)}}{\eta_1^2(c)}\eta_1^2(x)e^{-Q(x)}dx=2\frac{e^{Q(c)}}{\eta_1^2(c)}\eta_1^2(x)\mu(dx).
\end{equation}
We have proved that $\eta_1\in\mathbb{L}^2(\mu)$, i.e. $\int_{0}^{\infty}\eta_1^2(x)\mu(dx)<\infty$, which implies $\int_{0}^{\infty}\eta_1^{-2}(z)e^{Q(z)}dz=\infty$. Then $\Lambda^Y(\infty)=\infty$ and from (\ref{4.5}) we obtain $m(0,\infty)<\infty$.
\par
Let $T^{Y}_{a}:=\inf\{t>0:Y_t=a\}$ be the hitting time of $a$ for the process $Y$. For any $x,a\in(0,\infty)$, we know that the process $Y$ is positive recurrent when $\mathbb{E}_x(T^{Y}_{a})<\infty$. By using the formulas on page 353 in \cite{KS88}, we deduce $Y$ is positive recurrent. Therefore, $X^\tau$ is $R$-positive.
\qed
\end{pf}

\section{Examples}
\label{sect5}
In this section, we will illustrate our results with the following examples. Moreover, the second example is also given to exhibit the usefulness of the results.\\

{\bf{Example 1.}}
The first example we consider is the diffusion
$$dX_t=dB_t-adt,~~~~~~X_0=x>0,$$
where $a>0$ is constant. In this case, $q(x)=a$, $Q(y)=\int_{0}^{y}2adx=2ay$, $\Lambda(x)=\int_{0}^{x}e^{Q(y)}dy=\frac{1}{2a}(e^{2ax}-1)$. Then it is direct to check that
$\Lambda(\infty)=\infty, S=\int_{0}^{\infty}e^{Q(y)}\left(\int_y^{\infty}e^{-Q(z)}dz\right)dy=\int_{0}^{\infty}e^{2ay}\cdot\frac{1}{2a}e^{-2ay}dy=\infty$.
\par
Consider $\eta_\lambda$, the solution of
\begin{equation}
\label{5.1}
\frac{1}{2}\upsilon''(x)-a\upsilon'(x)=-\lambda\upsilon(x),~~~~\upsilon(0)=0,~~~\upsilon'(0)=1.
\end{equation}
If $a^2-2\lambda>0$, we have $\eta_\lambda(x)=\frac{1}{2\sqrt{a^2-2\lambda}}(e^{(a+\sqrt{a^2-2\lambda})x}-e^{(a-\sqrt{a^2-2\lambda})x})$,
then in this case, for any $x>0,~\eta_\lambda(x)>0$.
If $a^2-2\lambda=0$, we have $\eta_\lambda(x)=xe^{a x}$, then in this case, for any $x>0,~\eta_\lambda(x)>0$.
If $a^2-2\lambda<0$, we have $\eta_\lambda(x)=\frac{1}{\sqrt{2\lambda-a^2}}e^{a x}\sin(\sqrt{2\lambda-a^2}x)$
, then in this case, for any $x>0,~\eta_\lambda(x)$ has to change its sign.\par
Hence $\lambda_c=\frac{a^2}{2}$. By Proposition \ref{pro 3.1}, for any $0<\lambda\leq\lambda_c$,
$$d\nu_\lambda=2\lambda\eta_\lambda d\mu=2\lambda\eta_\lambda e^{-2ay}dy$$
is a QSD. In particular, the minimal QSD is $\nu_{\lambda_c}(dy)=a^2ye^{-ay}dy$. This result is in accordance with \cite{MM94}.\\

{\bf{Example 2.}}
The second example we consider is the Ornstein-Uhlenbeck process
$$dX_t=dB_t-aX_tdt,~~~~~~X_0=x>0,$$
where $a>0$ is constant. In this case, $q(x)=ax$, $Q(y)=\int_{0}^{y}2axdx=ay^2$, $\Lambda(x)=\int_{0}^{x}e^{Q(y)}dy=\int_{0}^{x}e^{ay^2}dy$. From this we have the following behaviors at $\infty$:
\begin{eqnarray*}
\int_{0}^{x}e^{ay^2}dy\underset{x\rightarrow\infty}{\sim}\frac{1}{2ax}e^{ax^2}~~~~~\mathrm{and}~~~~~\int_{x}^{\infty}e^{-ay^2}dy\underset{x\rightarrow\infty}{\sim}\frac{1}{2ax}e^{-ax^2}.
\end{eqnarray*}
Then, straightforward calculations show that
$$\Lambda(\infty)=\infty~~\mathrm{and}~~S=\int_{0}^{\infty}e^{Q(y)}\left(\int_y^{\infty}e^{-Q(z)}dz\right)dy=\infty.$$
\par
Consider $\eta_\lambda$, the solution of
\begin{equation}
\label{5.2}
\frac{1}{2}\upsilon''(x)-ax\upsilon'(x)=-\lambda\upsilon(x),~~~~\upsilon(0)=0,~~~\upsilon'(0)=1.
\end{equation}
\par For (\ref{5.2}), we know from \cite[Lemma 3.6]{LM00} that $\{\lambda~|~\varphi_\lambda(\cdot) ~\mathrm{does~ not~ change ~sign}\}=(-\infty,a]$ and for any $\lambda\in(0,a]$, $\int_0^{\infty}\varphi_\lambda(x)dx<\infty$. \par
Hence $\lambda_c=a$. By Proposition \ref{pro 3.1}, for any $0<\lambda\leq\lambda_c$,
$$d\nu_\lambda=2\lambda\eta_\lambda d\mu=2\lambda\eta_\lambda e^{-ay^2}dy$$
is a QSD. This result is in accordance with \cite{LM00}.
\par
By using the above asymptotic relation, it is direct to check that
\begin{eqnarray*}
\lim\limits_{x\rightarrow\infty}\mu([x,\infty))\int_{0}^{x}e^{Q(y)}dy=\lim\limits_{x\rightarrow\infty}\frac{1}{4a^2x^2}=0.
\end{eqnarray*}
Therefore, by Theorem \ref{thm4.1} it follows that the Ornstein-Uhlenbeck process killed at 0 is $R$-positive.

\section*{Acknowledgements}
The authors would like to thank the referee for helpful comments and for pointing out the work of Pinsky \cite{PRG09}.
The first author also would like to thank the Centro de Modelamiento Matem\'{a}tico of Universidad de Chile, where part of this work was done, for their kind hospitality. The work is supported by the National Natural Science Foundation of China (Grant No.11371301) and the Science and Technology Planning Project of Hunan Province of China (Grant No. 2012FJ4093).



%% References


%\bibliographystyle{elsarticle-num}

\begin{thebibliography}{00}


\bibitem{CCLMMS09}
P. Cattiaux, P. Collet, A. Lambert, S. Mart\'{\i}nez, S. M\'{e}l\'{e}ard, J. San Mart\'{\i}n,
Quasi-stationary distributions and diffusion models in population dynamics,
Ann. Probab. 37 (2009) 1926--1969.%

\bibitem{MFC00}
M.F. Chen,
Explicit bounds of the first eigenvalue,
Sci. in China, Ser. A. 43 (2000) 1051--1059.%

\bibitem{CL55}
E.A. Coddington, N. Levinson,
Theory of Ordinary Differential Equtions, McGraw-Hill, New York, 1955.%

\bibitem{CMS95}
P. Collet, S. Mart\'{\i}nez, J. San Mart\'{\i}n,
Asymptotic laws for one-dimensional diffusions conditioned to nonabsorption,
Ann. Probab. 23 (1995) 1300--1314.%

\bibitem{IW89}
N. Ikeda, S. Watanabe, Stochastic Differential Equations and Diffusion
Processes (North-Holland Mathematical Library {\bf24}), 2nd edn, North-Holland, Amsterdam, 1989.%

\bibitem{KS88}
I. Karatzas, S. Shreve,
Brownian Motion and Stochastic Calculus, Springer, New York, 1988.%

\bibitem{KT81}
S. Karlin, H.M. Taylor, A Second Course in Stochastic Processes, 2nd edn, Academic Press, New York, 1981.%

\bibitem{KS12}
M. Kolb, D. Steinsaltz, Quasilimiting behavior for one-dimensional diffusions with killing, Ann. Probab. 40 (2012) 162--212.%

\bibitem{LC12}
J. Littin,
Uniqueness of quasistationary distributions and discrete spectra when $\infty$ is an entrance boundary and 0 is singular,
J. Appl. Probab. 49 (2012) 719--730.%

\bibitem{LM00}
M. Lladser, J. San Mart\'{\i}n, Domain of attraction of the quasi-stationary distributions for the Ornstein-Uhlenbeck process, J. Appl. Probab. 37 (2000) 511--520.%

\bibitem{M61}
P. Mandl, Spectral theory of semi-groups connected with diffusion processes and its application, Czech. Math. J. 11 (1961) 558--569.%

\bibitem{MM94}
S. Mart\'{\i}nez, J. San Mart\'{\i}n, Quasi-stationary distributions for a Brownian motion with drift and associated limit laws, J. Appl. Probab. 31 (1994) 911--920.%

\bibitem{MPM98}
S. Mart\'{\i}nez, P. Picco, J. San Mart\'{\i}n, Domain of attraction of quasi-stationary distributions for the Brownian motion with drift, Adv. Appl. Probab. 30 (1998) 385--408.%

\bibitem{MM01}
S. Mart\'{\i}nez, J. San Mart\'{\i}n, Rates of decay and $h$-processes for one dimensional diffusions conditional on non-absorption, J. Theoret. Probab. 14 (2001) 199--212.%

\bibitem{MM04}
S. Mart\'{\i}nez, J. San Mart\'{\i}n, Classification of killed one-dimensional diffusions, Ann. Probab. 32 (2004) 530--552.%

\bibitem{PRG09}
R.G. Pinsky, Explicit and almost explicit spectral calculations for diffusion operators, J. Funct. Anal. 256 (2009) 3279--3312.%

\bibitem{SE07}
D. Steinsaltz, S.N. Evans, Quasistationary distributions for one-dimensional diffusions with killing, Trans. Amer. Math. Soc. 359 (2007) 1285--1324 (electronic).%

\bibitem{TT79}
P. Tuominen, R.L. Tweedie, Exponential decay and ergodicity of general Markov processes and their discrete skeletons, Adv. Appl. Probab. 11 (1979) 784--803.%










\end{thebibliography}






%\end{thebibliography}

\end{document}